# Scaling Qualitative Probability

Mark Burgin

University of California, Los Angeles
405 Hilgard Ave.
Los Angeles, CA 90095

**Abstract**

There are different approaches to qualitative probability, which includes subjective probability. We developed a representation of qualitative probability based on relational systems, which allows modeling uncertainty by probability structures and is more coherent than existing approaches. This setting makes it possible proving that any comparative probability is induced by some probability structure (Theorem 2.1), that classical probability is a probability structure (Theorem 2.2) and that inflated, i.e., larger than 1, probability is also a probability structure (Theorem 2.3). In addition, we study representation of probability structures by classical probability.

## 1. Introduction

Qualitative probability and its particular case - subjective probability – plays important role in psychology, decision making and cognitive sciences. However, the mainstream of probability theory and statistics treats probability as a numeric function or measure with specific properties (cf., for example, (Devore, 2015; Hogg, et al, 2014; Walpole, et al, 2011)). Kolmogorov rigorously represented these properties by a system of axioms treating probability as an additive measure on a system of random events (Kolmogorov, 1933). Later researchers suggested other systems of axioms for the conventional probability theory (cf., for example, (Whittle, 2000)), as well as for extended probability theory, which included negative probabilities (Burgin, 2009; 2016), and for larger than 1 (inflated) probabilities (Burgin and Meissner, 2012).

At the same time, to reflect uncertainty and probabilistic aspects of decision-making and other cognitive processes, various researchers developed and applied the theory of *qualitative probability* as an approach to mathematical modeling and investigation of uncertainty. While many properties of classical probability are consequences of properties of real numbers, in which classical probabilities take values, qualitative probability demands a more sophisticated description. As a result, different axiomatic systems have been suggested for qualitative probabilities (cf., for example, (Bernstein, 1917; Ramsey, 1931; de Finetti, 1931; 1937; Koopman, 1940; Cox, 1946; Savage, 1954; Fishburn, 1986; Watson and Buede, 1987)).

According to the general approach to qualitative probability, it is represented by comparative probability relations on sets of events making possible comparison of probability for different events. Usually, three probability relations are used:

- The strict comparative probability relation $\succ$, for which the expression $A \succ B$ means "the event $A$ is more probable than the event $B$".

- The non-strict comparative probability relation $\succcurlyeq$, for which the expression $A \succcurlyeq B$ means "the event $A$ is more probable than the event $B$ or has the same probability".

- The equivalence probability relation $\asymp$, for which the expression $A \asymp B$ means either comparative equiprobability or absence of significant difference between probabilities of events $A$ and $B$.

However, these relations are not between events but between probabilities of these events. To reflect this peculiarity, we build a probability theory that distinguishes events and their probabilities in this paper. In Section 2, probability structures are introduced as an extension and generalization of probability functions and their properties are studied. Section 3 explores representation of some probability structures by other probability structures, for example, representation of inflated probability by classical probability or some classes of total probability structures by classical probability. The structural approach to probability allows achieving unification of different forms, types and kinds of probability.

## 2. Probability structures

Let us consider a set $U$, which is called the *universe of events* and elements of which are called *elementary events*. Subsets of the set $U$ are called *events*. In what follows, we assume that the set $U$ is finite, i.e., $U = \{w_1, w_2, w_3, \ldots, w_n\}$. Probability structures with infinite universes of events are studied elsewhere.

Usually, probability is defined not for all events but only for random events. However, the concept of randomness has different interpretations and formalizations. That is why to include the classical approach as well as other approaches to probability, we specify a set **F** of *permissible*, e.g., random, events, which consists of some subsets of *U* and is defined by its formal properties. Elements from the set $2^U$ of all subsets of *U* that belong to **F** are called *elementary permissible events*.

This approach allows treating different types of random events. Even more, it makes possible treating all events as permissible as it is done in some theories of subjective probability (cf., (Fishburn, 1986)) and in the theory of hyperprobability (Burgin and Krinik, 2009).

In some cases, instead of probability of events, probability of system states is studied (cf., for example, (Fishburn, 1986)). However, it is possible to treat the probability of a state *q* of a system *R* as the probability of the event such as finding the system *R* in the state *q* and the probability of a set {$q_1$, $q_2$, $q_3$, …, $q_k$} of states of the system *R* as the probability of the event such as finding the system *R* in the states $q_1$, $q_2$, $q_3$, …, $q_k$.

Sometimes instead of probability of events, probability of statements, e.g., propositions or predicates, is studied. In this case, it is possible to treat the probability of a statement *P* as the probability of the event of finding that the statement *P* is true or simply, the probability of the event such as finding that the statement *P* is true.

For the development of probability theory, it is important that the set **F** has a flexible algebraic structure.

**Axiom SIP1** (*Boolean structure*). **F** is a Boolean algebra with *U* as its element.

A more advanced structures are set algebras and set fields. We remind that a collection of sets is a *set algebra* if it is closed with respect to union, intersection and difference of sets and a set algebra **B** closed with respect to complements of its elements is called a *set field* (Kolmogorov and Fomin, 1999). For any set algebra **A**, the empty set Ø belongs to **A** and for any set field **B** in *U*, the set *U* belongs to **B**.

**Axiom SIP2** (*Algebraic structure*). **F** is a set algebra that has *U* as its element.

**Axiom SIP3** (*Field structure*). **F** is a set field in *U*.

In some approaches to subjective probability, pseudocomplemented distributive lattices or orthocomplemented lattices are also considered (Narens, 2007; 2014).

To define probability of permissible events, we consider a mathematical structure $X$, which is called a *probability scale*. For comparative probabilities, the scale $X$ is a partially ordered set with the strict order relation $\succ$, for which the expression $a \succ b$ means "the element $a$ is larger than the element $b$".

The relation $\succ$ induces the non-strict order relation $\succcurlyeq$, for which the expression $a \succcurlyeq b$ means "the element $a$ is larger than the element $b$" or these elements coincide, i.e., $a = b$.

**Examples:**
1. The scale of the conventional probability is the interval [0, 1] of real numbers with the natural ordering of real numbers (Kolmogorov, 1933).
2. The scale of extended probability is the interval [-1, 1] of real numbers with the natural ordering of real numbers (Bartlett, 1945; Burgin, 2009).
3. The scale of inflated probability is the interval [0, $k$] where $k$ is a positive real number larger than 1 with the natural ordering of real numbers (Burgin and Meissner, 2012).

Note that there are different concepts of a scale in other areas.

**Definition 2.1.** A mapping $P: \mathbf{F} \to X$ is called a *positive probability structure*, if $\mathbf{F}$ satisfies axioms SIP1, SIP2 or SIP3 and $P$ satisfies the following axioms:

**Positive Non-triviality Axiom (PNA):** If $U \neq \emptyset$, then $P(U) \succ P(\emptyset)$.

**Inclusion Correlation Axiom (ICA):** For any permissible events $A$ and $B$, if $A \supseteq B$, then $P(A) \succcurlyeq P(B)$.

**Positive Additivity Axiom (PAA):** For any permissible events $A$, $B$ and $C$, if $A \cap C = B \cap C = \emptyset$, then
$$P(A) \succ P(B) \Leftrightarrow P(A \cup C) \succ P(B \cup C)$$

We see that a positive probability structure actually consists of two mathematical structures, i.e., sets with relations, and a mapping of the first structure $\mathbf{F}$ into the second structure $X$. It means that a positive probability structure forms a structure of the second level and this structure is called a named set (Burgin, 2011).

**Examples:**
1. The conventional probability is a mapping of set algebra $\mathbf{F}$ into the interval [0, 1] of real numbers with the natural ordering of real numbers (Kolmogorov, 1933), i.e., it is a probability structure as it is proved in Theorem 2.1.

2. The extended probability is a mapping of set algebra **F** into the interval [-1, 1] of real numbers with the natural ordering of real numbers (Bartlett, 1945; Burgin, 2009), i.e., it is also a probability structure, which is not positive. Such structures are studied elsewhere.
3. The scale of inflated probability is a mapping of set algebra **F** into the interval [0, *k*] where *k* is a positive real number larger than 1 with the natural ordering of real numbers (Burgin and Meissner, 2012), i.e., it is a probability structure as it is proved in Theorem 2.2.
4. A probability measure is a mapping of set algebra **F** into the interval [0, 1] of real numbers with the natural ordering of real numbers (cf., for example, (Davidson and Suppes, 1956; Schmeidler 1984)), i.e., it is a probability structure.

One more example of probability structures is orthoprobability functions studied by Narens (2007; 2014).

Let $\mathbf{L} = \{L; \cap, \cup, \perp, U, \emptyset \}$ is an orthocomplemented lattice of subsets of the set *U*.

**Definition 2.2.** A mapping $P: \mathbf{L} \to X$ mapping is called an *orthoprobability function* (Narens, 1980; 2014) if it satisfies the following axioms:

**NA 1.** $P(\mathbf{L}) = 1$ and $P(\emptyset) = 0$.

**NA 2.** For any *A* and *B* from **L**, if $B \cap A = \emptyset$, then

$$P(A \cup B) = P(A) + P(B)$$

Let us study properties of positive probability structures.

**Proposition 2.1.** For any permissible event *A*, the inequality $A \neq U$ implies $P(U) \succcurlyeq P(A)$.

Indeed, if $A \subseteq U$, then $U = A \cup B$ and $B \cap A = \emptyset$. By Axiom ICA, we have $P(B) \succcurlyeq P(\emptyset)$. Then by Axiom PPA, we have $P(U) = P(A \cup B) \succcurlyeq P(\emptyset \cup A) = P(A)$, i.e., $P(U) \succcurlyeq P(A)$.

Positive probability structures also satisfy additional axioms utilized in the theory of qualitative probability.

**Proposition 2.2.** For an arbitrary positive probability structure *P*, the following axioms are true:

**Asymmetry Axiom (ASA):** For any permissible events *A* and *B*, if $P(A) \succ P(B)$, then it is not $P(B) \succ P(A)$.

**Inclusion Monotonicity Axiom (IMA)**: For any permissible events $A$, $B$ and $C$, if $A \supseteq B$ and $P(B) \succcurlyeq P(C)$, or $B \supseteq C$ and $P(A) \succcurlyeq P(B)$, then $P(A) \succcurlyeq P(C)$.

*Proof.* 1. Axiom ASA: Indeed, a strict partial order is an asymmetric relation.

2. Axiom IMA: Indeed, by Axiom MA, which is already proved, $A \supseteq B$ implies $A \succcurlyeq B$ and $A \succcurlyeq B$ and $B \succcurlyeq C$ imply $A \succcurlyeq C$ because by Lemma 2.1, $\succcurlyeq$ is a partial order.

Proposition is proved.

**Lemma 2.1.** For any permissible events $A$ and $B$, we have $P(A \cup B) \succcurlyeq P(A)$ and $P(A) \succcurlyeq P(A \cap B)$.

*Proof.* 1. By the laws of set theory, $A \cup B = A \cup (B \setminus A)$ and $(B \setminus A) \cap A = \emptyset$. By Axiom ICA, we have $P(B \setminus A) \succcurlyeq P(\emptyset)$. Then by Axiom PPA, we have $P(A \cup B) = P(A \cup (B \setminus A)) \succcurlyeq P(A \cup \emptyset) = P(A)$, i.e., $P(A \cup B) \succcurlyeq P(A)$.

2. By the laws of set theory, $A = (A \cap B) \cup (A \setminus B)$ and $(A \setminus B) \cap (A \cap B) = \emptyset$. By Axiom ICA, we have $P(A \setminus B) \succcurlyeq P(\emptyset)$. Then by Axiom PPA, we have $P(A) = P((A \cap B) \cup (A \setminus B)) \succcurlyeq P((A \cap B) \cup \emptyset) = P(A \cap B)$, i.e., $P(A) \succcurlyeq P(A \cap B)$.

Lemma is proved.

Taking a positive probability structure $P$, it is possible to compare probabilities of the events from **F**. Namely, for events $A$ and $B$, the expression $P(A) \succ P(B)$ means "the event $A$ is more probable than the event $B$", the expression $P(A) \succcurlyeq P(B)$ means "the event $A$ is more probable than the event $B$ or has the same probability" and the expression $P(A) = P(B)$ means either absence of significant difference between or equality of probabilities of events $A$ and $B$.

Besides, relations in the set $X$ induce corresponding relations in the set **F**:

$$A \succ B \text{ if } P(A) \succ P(B)$$

$$A \succcurlyeq B \text{ if } P(A) \succcurlyeq P(B)$$

In addition, we can define $A \asymp B$ when $P(A) = P(B)$.

This allows us to compare positive probability structures and comparative probabilities studied in the literature.

The most popular system of qualitative probability axioms was introduced by de Finetti and used by many other authors (cf., for example, (Savage (1954; Kaplan and Fine, 1977; Wakker, 1981; French, 1982; Fishburn, 1983; 1986; Suppes, 2014)). This system consists of the following axioms.

**Weak Order Axiom (WOA):** The relation $\succ$ is asymmetric and transitive.

**Non-triviality Axiom (NTA)**: If $U \neq \emptyset$, then $U \succ \emptyset$.

**Non-negativity Axiom (NNA)**: For any event $A$, $A \succcurlyeq \emptyset$.

**Additivity Axiom (AA)**: For any permissible events $A$, $B$ and $C$, if $A \cap C = B \cap C = \emptyset$, then

$$A \succ B \Leftrightarrow A \cup C \succ B \cup C$$

Let us consider some consequences of these axioms.

**Lemma 2.3.** The relation $\succcurlyeq$ is a partial order.

Axiom WOA directly implies the following axioms.

**Strict Transitivity Axiom (STA)**: For any permissible events $A$ and $B$, if $A \succ B$ and $B \succ C$, then $A \succ C$.

**Non-strict Transitivity Axiom (NTA)**: For any permissible events $A$, $B$ and $C$, if $A \succcurlyeq B$ and $B \succcurlyeq C$, then $A \succcurlyeq C$.

**Lemma 2.4.** Axiom WOA implies the following axiom used for qualitative probability:

**Asymmetry Axiom (ASA)**: For any permissible events $A$ and $B$, if $A \succ B$, then it is not $B \succ A$.

Indeed, a strict partial order is an asymmetric relation.

**Lemma 2.4.** Axioms NNA and AA imply the following axioms used for qualitative probability:

**Monotonicity Axiom (MA)**: For any permissible events $A$ and $B$, if $A \supseteq B$, then $A \succcurlyeq B$.

**Inclusion Monotonicity Axiom (IMA)**: For any permissible events $A$, $B$ and $C$, if $A \supseteq B$ and $B \succcurlyeq C$, or $B \supseteq C$ and $A \succcurlyeq B$, then $A \succcurlyeq C$.

*Proof.* 1. Axiom MA: If $A \supseteq B$, then $A = B \cup C$ and $B \cap C = \emptyset$. By Axiom NNA, we have $C \succcurlyeq \emptyset$. Then by Lemma 2.2, we have $A = C \cup B \succcurlyeq \emptyset \cup B = B$, i.e., $A \succcurlyeq B$. This gives us Axiom MA.

2. Axiom IMA: Indeed, by Axiom MA, which is already proved, $A \supseteq B$ implies $A \succcurlyeq B$ and $A \succcurlyeq B$ and $B \succcurlyeq C$ imply $A \succcurlyeq C$ because by Lemma 2.1, $\succcurlyeq$ is a partial order. This gives us Axiom IMA.

Lemma is proved.

Let us consider two positive probability structures $P: \mathbf{F} \to X$ and $Q: \mathbf{H} \to X$.

**Definition 2.3.** a) A positive probability structure $P$ is a *substructure* of a positive probability structure $Q$ if $\mathbf{F}$ is a subset of $\mathbf{H}$ and for any elements $e$ and $h$ from $\mathbf{F}$, $P(e) \succ P(h)$ if and only if $Q(e) \succ Q(h)$.

b) A positive probability structure $P$ is a *robust substructure* of a positive probability structure $Q$ if the following diagram in which $m: \mathbf{F} \to \mathbf{H}$ is an inclusion is commutative, i.e., for any event $e$ from $\mathbf{F}$, we have $Q(m(e)) = P(e)$,

$$\begin{array}{c} Q \\ \mathbf{H} \to X \\ m \nwarrow \quad \nearrow P \\ \mathbf{F} \end{array} \qquad (1)$$

**Lemma 2.5.** If a positive probability structure $P$ is a substructure of a positive probability structure $Q$ and a positive probability structure $Q$ is a substructure of a positive probability structure $R$, then a positive probability structure $P$ is a substructure of a positive probability structure $R$.

*Proof.* Let us consider three positive probability structures $P: \mathbf{F} \to X$, $P: \mathbf{G} \to X$ and $Q: \mathbf{H} \to X$ such that the structure $P$ is a substructure of the structure $Q$ and the structure $Q$ is a substructure of the structure $R$. It means that $\mathbf{F}$ is a subset of $\mathbf{H}$ and $\mathbf{H}$ is a subset of $\mathbf{G}$. By properties of sets, $\mathbf{F}$ is a subset of $\mathbf{G}$. Besides, for any elements $e$ and $h$ from $\mathbf{F}$, $P(e) \succ P(h)$ if and only if $Q(e) \succ Q(h)$ and $Q(e) \succ Q(h)$ if and only if $R(e) \succ R(h)$ because elements $e$ and $h$ also belong to **H**. Consequently, $P(e) \succ P(h)$ if and only if $R(e) \succ R(h)$. It means that a positive probability structure $P$ is a substructure of a positive probability structure $R$.

Lemma is proved.

**Lemma 2.6.** If a positive probability structure $P$ is a robust substructure of a positive probability structure $Q$ and a positive probability structure $Q$ is a robust substructure of a positive probability structure $R$, then a positive probability structure $P$ is a robust substructure of a positive probability structure $R$.

Indeed, combining the commutative diagram (1) with the commutative diagram (2), we obtain the commutative diagram (3), in which $l$ is an inclusion.

$$\begin{array}{c} R \\ \mathbf{G} \to X \\ k \nwarrow \nearrow Q \\ \mathbf{H} \end{array} \qquad (2)$$

$$\begin{array}{c} R \\ \mathbf{G} \to X \\ l \nwarrow \nearrow P \\ \mathbf{F} \end{array} \qquad (3)$$

Definitions imply the following result.

**Lemma 2.7.** Any robust substructure of a positive probability structure is a substructure of the same probability structure.

Let us consider some useful classes of positive probability structures.

**Definition 2.4.** A positive probability structure $P$: $\mathbf{F} \to X$ is called *rigid*, if it satisfies the following axiom:

**Equality Correlation Axiom (ECA)**: For any permissible events $A$, $B$ and $C$, if $A \cap C = B \cap C = \emptyset$, then

$$P(A) = P(B) \Leftrightarrow P(A \cup C) = P(B \cup C)$$

Definitions imply the following result.

**Lemma 2.8.** A (robust) substructure of a rigid positive probability structure is a rigid positive probability structure.

**Proposition 2.2.** For any permissible events A, B and C in a rigid positive probability structure P, if $A \cap C = B \cap C = \emptyset$, then

$$P(A) \succcurlyeq P(B) \Leftrightarrow P(A \cup C) \succcurlyeq P(B \cup C).$$

Indeed, $P(A) \succcurlyeq P(B)$ means that either $P(A) = P(B)$ or $P(A) \succ P(B)$. In the first case, $P(A) = P(B) \Leftrightarrow P(A \cup C) = P(B \cup C)$ by Axiom ECA. In the second case, $P(A) \succ P(B) \Leftrightarrow P(A \cup C) \succ P(B \cup C)$ by Axiom PAA.

**Definition 2.5.** A positive probability structure $P: \mathbf{F} \to X$ is called *complete*, if $\succ$ is a total strict order in the image of the set $\mathbf{F}$.

**Proposition 2.3.** A positive probability structure is complete if and only if any of its substructures is complete.

Indeed, any positive probability structure has the same scale as any of its substructures. Therefore, if a substructures is complete, then the whole probability structure is also complete.

Definitions imply the following result.

**Proposition 2.4.** In a complete positive probability structure, the following axiom is true

**Positive Completeness Axiom (PCTA):** For any permissible event A, either $P(A) \succcurlyeq P(B)$ or $P(B) \succcurlyeq P(A)$.

Let us consider some classes of positive probability structures.

**Definition 2.6.** A positive probability structure $P: \mathbf{F} \to X$ is called *total*, if $\mathbf{F}$ is the power-set $2^U$ of U, i.e., the set of all subsets of U.

**Proposition 2.5.** A complete positive probability structure that satisfies Axiom SIP1 is total if and only if $\mathbf{F}$ contains all one-element subsets of U.

Indeed, any subset of U is the union of one-element subsets of U.

Assuming that the relation $\asymp$ is an equivalence relation correlated with the relations defined on events from $\mathbf{F}$ by comparative probability, we can prove the following result.

**Theorem 2.1.** a) Any acyclic comparative probability relation is induced by some positive probability structure P.

b) If a comparative probability relation $\succcurlyeq$ is induced by a complete positive probability structure P, the following axiom is true:

**Completeness Axiom (CTA)**: For any event $A$, $A \succcurlyeq B$ or $B \succcurlyeq A$.

To prove this theorem, we consider the set **F** of permissible events with the comparative probability relation $\succ$ and equivalence probability relation $\asymp$ in it and take the factor-set $\mathbf{F}/\!\asymp$ with the induced by $\asymp$ partial order as the scale $X$ of the positive probability structure $P$. Then this order induces relations $\succ$ and $\asymp$, while Axioms WOA, NTA, NNA and AA imply Axioms PNA, ICA and PAA for the positive probability structure $P$.

Totality of the order in $X$ implies Axiom CTA.

We call a probability function (distribution) $P: \mathbf{F} \to [0, 1]$ classical if it satisfies Kolmogorov axioms (Kolmogorov, 1933). When $\mathbf{F} \subseteq 2^\Omega$, they have the form:

**K 1.** (*Non-negativity*) $P(A) \geq 0$ for all $A \in \mathbf{F}$.

**K 2.** (*Normalization*) $P(\Omega) = 1$.

**K 3.** (*Finite additivity*)

$$P(A \cup B) = P(A) + P(B)$$

for all sets $A, B \in \mathbf{F}$ such that

$$A \cap B = \emptyset.$$

There is an intrinsic relation between classical probability functions (distributions) and positive probability structures.

**Theorem 2.2.** Any classical probability function (distribution) is a complete rigid positive probability structure.

*Proof.* To prove this result, we need to show that a classical probability function $P$ satisfies only axioms PNA, ICA, PAA and ECA because any classical probability function is defined on a set algebra and the natural order of real numbers is linear (total).

Positive Non-triviality Axiom (PNA) directly follows from Axiom K2.

Inclusion Correlation Axiom (ICA) follows from Axioms K1 and K3. Indeed, for any events $A$ and $B$, if $A \supseteq B$, then $A = B \cup C$ where $C = A \setminus B$. Thus, $B \cap C = \emptyset$. Then by Axiom K3, we have $P(A \cup B) = P(C) + P(B)$ and $P(A) \geq P(B)$ as $P(C) \geq 0$ by Axiom K1. This gives us Axiom ICA.

Positive Additivity Axiom (PAA) directly follows from Axiom K3. Indeed, for any permissible events $A, B$ and $C$, if $A \cap C = B \cap C = \emptyset$, we have

$$P(A \cup C) = P(A) + P(C)$$

and

$$P(B \cup C) = P(B) + P(C)$$

Then by properties of order relations, we have

$$P(A) > P(B) \Leftrightarrow P(A) + P(C) > P(B) + P(C) \Leftrightarrow P(A \cup C) > P(B \cup C)$$

This gives us Axiom PAA.

Equality Correlation Axiom (ECA) is also true. Indeed, for any permissible events *A*, *B* and *C*, by Axiom FA, if $A \cap C = B \cap C = \emptyset$, then

$$P(A \cup C) = P(A) + P(C)$$

and

$$P(B \cup C) = P(B) + P(C)$$

Thus,

$$P(A) = P(B) \Leftrightarrow P(A \cup C) = P(B \cup C)$$

Theorem is proved.

**Corollary 2.1.** The classical probability theory is a subtheory of qualitative probability theory in the form of probability structures.

In (Burgin and Meissner, 2012), inflated, i.e., larger than 1, probabilities were rigorously introduced in a mathematically consistent way, studied and applied to modeling financial processes. Note that larger than 1 probabilities were used in physics without mathematical grounding by such physicists as Dirac, Weisskopf and Wigner.

**Theorem 2.3.** Any inflated probability function (distribution) is a complete rigid positive probability structure.

*Proof.* To prove this result, we need to show that an inflated probability function *P* with the set Ω of elementary events satisfies axioms PNA, ICA, PAA and ECA because the natural order of real numbers is linear (total). Inflated probability functions are characterized by two groups of axioms – physical axioms and mathematical axioms.

Physical axioms are:

**Decomposability Axiom.** An event *A* occurs in a trial if and only if some elementary event from *A* occurs in this trial.

**Boundary Axiom.** The number of occurrences of an elementary event $w_i$ from *A* in one trial cannot be larger than the multiplicity $n_i$ of this event in Ω.

**Occurrence Axiom.** In any trial, at least one elementary event occurs.

**Finite Occurrence Axiom.** In any sequence of trials, there is only a finite number of trials when at least one elementary event does not occur.

However, we also need a weaker condition.

**Infinite Occurrence Axiom (IOA).** In any sequence of trials, there are an infinitely many trials when at least one elementary event occurs.

**Finite Occurrence Condition (FOC).** In any sequence of trials, there only a finite number of trials when an elementary event from *A* occurs.

**Outcome Axiom** (**OA**). In each trial, exactly one elementary event happens.

**Single Type Axiom (STA).** In each trial, elementary events of only one type happen and their number cannot be larger than the multiplicity of this element in the multiset Ω.

Mathematical axioms are:

**Axiom AS** (*Algebraic structure*). **F** is a set algebra that has Ω as its element.

**Axiom UP** (*Upper normalization*). $0 \leq P(\Omega) \leq m$ for some natural number *m*.

**Axiom LP** (*Lower normalization*). $P(\emptyset) = 0$ for some natural number *m*.

**Axiom FA** (*Finite additivity*)

$$P(A \cup B) = P(A) + P(B)$$

for all sets $A, B \in \mathbf{F}$ such that

$$A \cap B \equiv \emptyset$$

Note that Axiom LP follows from Outcome Axiom OA considered in (Burgin and Meissner, 2012).

As probability structures are characterized only by mathematical axioms, we need only mathematical axioms of inflated probabilities to prove this theorem. Thus, we check axioms of probability structures.

Any interval [0, *m*] of real numbers, which is the scale of an inflated probability function, is linearly (totally) ordered. Therefore, the basic condition of probability structures – partial order in the scale – is true.

Let us check other axioms.

Positive Non-triviality Axiom (PNA) directly follows from Axioms OA and UP.

Inclusion Correlation Axiom (ICA) follows from Axiom FA. Indeed, for any events $A$ and $B$, if $A \supseteq B$, then $A = B \cup C$ where $C = A \setminus B$. Thus, $B \cap C = \emptyset$. Then by Axiom FA, we have $P(A \cup B) = P(C) + P(B)$ and $P(A) \geq P(B)$ as $P(C) \geq 0$ by Axiom UP. This gives us Axiom ICA.

Positive Additivity Axiom (PAA) directly follows from Axioms OA and UP. Indeed, for any permissible events $A$, $B$ and $C$, if $A \cap C = B \cap C = \emptyset$, we have

$$P(A \cup C) = P(A) + P(C)$$

and

$$P(B \cup C) = P(B) + P(C)$$

Then by properties of order relations, we have

$$P(A) > P(B) \Leftrightarrow P(A) + P(C) > P(B) + P(C) \Leftrightarrow P(A \cup C) > P(B \cup C)$$

This gives us Axiom PAA.

Equality Correlation Axiom (ECA) is also true. Indeed, for any permissible events $A$, $B$ and $C$, by Axiom FA, if $A \cap C = B \cap C = \emptyset$, then

$$P(A \cup C) = P(A) + P(C)$$

and

$$P(B \cup C) = P(B) + P(C)$$

Thus,

$$P(A) = P(B) \Leftrightarrow P(A \cup C) = P(B \cup C)$$

Theorem is proved.

## 3. Representation of probability structures

Representation is an important relation between probability structures.

Let $X$ and $Y$ by partially ordered sets with the order relation $<$. We remind the following concept.

**Definition 3.1.** A mapping $f: X \to Y$ is called an *order homomorphism*, or simply, a *homomorphism*, if for any elements $a$ and $b$ from $X$, we have

$$a < b \text{ implies } f(a) < f(b)$$

Inclusion relation $\subseteq$ of sets is a partial order in systems of sets. Thus, Axiom ICA implies the following result.

**Proposition 3.1.** Any positive probability structure $P$ is an order homomorphism.

Homomorphisms or mappings that preserve definite structures play an important role in many mathematical theories.

**Definition 3.2.** A homomorphism $f: X \to Y$ is called a *monomorphism* if it is a one-to-one mapping into $Y$.

For instance, there is a monomorphism of whole numbers into rational numbers.

Let us consider two positive probability structures $P: \mathbf{F} \to X$ and $Q: \mathbf{F} \to Y$.

**Definition 3.3.** A positive probability structure $P$ is *represented* by a positive probability structure $Q$ if there is a homomorphism $f: X \to Y$ such that the following diagram is commutative, i.e., for any event $e$ from $\mathbf{F}$, we have $f(P(e)) = Q(e)$,

$$\begin{array}{c} P \\ \mathbf{F} \to X \\ {}_Q\searrow \swarrow_f \\ Y \end{array} \quad (4)$$

As a composition of homomorphisms is a homomorphism, we have the following result.

**Proposition 3.2.** If a positive probability structure $P: \mathbf{F} \to X$ is represented by a positive probability structure $Q: \mathbf{F} \to Y$ and the positive probability structure $Q$ is represented by a positive probability structure $T: \mathbf{F} \to Z$, then the positive probability structure $P$ is represented by the positive probability structure $T$.

Indeed, combining the commutative diagram (4) with the commutative diagram (5), we obtain the commutative diagram (6), in which $h$ is a homomorphism because the sequential composition of homomorphisms is a homomorphism.

$$\begin{array}{c} Q \\ \mathbf{F} \to Y \\ {}_T\searrow \swarrow_k \\ Z \end{array} \quad (5)$$

$$\begin{array}{c} P \\ \mathbf{F} \to Y \\ {}_T\searrow \swarrow_{h = k \circ f} \\ Z \end{array} \quad (6)$$

**Lemma 3.1.** If a positive probability structure $P: \mathbf{F} \to X$ is represented by a positive probability structure $Q: \mathbf{F} \to Y$, then any substructure $M: \mathbf{H} \to X$ of $P$ is represented by some substructure $N: \mathbf{H} \to Y$ of $Q$.

Indeed, the restriction of the mapping $Q$ on the subset $\mathbf{H}$ of the set $\mathbf{F}$ will give us the necessary representation.

Definitions imply the following result.

**Lemma 3.2.** If a positive probability structure $P: \mathbf{F} \to X$ is represented by a positive probability structure $Q: \mathbf{F} \to Y$, then $Q$ is total if and only if $P$ is total.

As in the Diagram (4) the image of $Q$ coincides with the image of $P \circ f$ and $f$ is an order homomorphism, we have the following result.

**Lemma 3.3.** If a complete positive probability structure $P: \mathbf{F} \to X$ is represented by a positive probability structure $Q: \mathbf{F} \to Y$, then $Q$ is also complete.

There different classes of probability structure representations.

**Definition 3.4.** A positive probability structure $P$ is *faithfully represented* by a positive probability structure $Q$ if there is commutative diagram (1) in which $f$ is a monomorphism.

As a composition of monomorphisms is a monomorphism, Proposition 3.1 implies the following result.

**Proposition 3.2.** If a positive probability structure $P$ is faithfully represented by a positive probability structure $Q$ and a positive probability structure $Q$ is faithfully represented by a positive probability structure $T$, then positive probability structure $P$ is faithfully represented by a positive probability structure $T$.

Rigidity is preserved by inverse faithful representations.

**Lemma 3.4.** If a positive probability structure $P: \mathbf{F} \to X$ is faithfully represented by a rigid positive probability structure $Q: \mathbf{F} \to Y$, then $P$ is also rigid.

Indeed, any equality of the form $Q(A) = Q(B)$ implies the equality $P(A) = P(B)$ by the definition of a faithful representation. Thus, validity of Axiom ECA for $Q$ implies its validity for $P$.

**Theorem 3.1.** It is possible to faithfully represent any inflated probability function by a classical probability function.

*Proof.* Let us consider an inflated probability function $P: \mathbf{F} \to [0, m]$ with the set $\Omega$ of elementary events. By Axiom UP, $0 \leq P(\Omega) \leq m$ for some natural number $m$. Then it is possible to build a mapping $Q: \mathbf{F} \to [0, 1]$ by the following rule: for any $A \in \mathbf{F}$, we define $Q(A) = P(A)/k$ where $k$ is equal to $P(\Omega)$. As for any $A \in \mathbf{F}$, we have $P(A) \leq P(\Omega)$ by Axiom FA, $Q$ maps $\mathbf{F}$ into the interval $[0, 1]$. By construction, we also have $P(A) \geq 0$ for all $A \in \mathrm{F}$ (Axiom K 1) and $P(\Omega) = 1$ (Axiom K 2). In addition, by Axiom FA, we have

$$Q(A \cup B) = P(A \cup B)/k = (P(A) + P(B))/k = P(A)/k + P(B)/k = Q(A) + Q(B)$$

i.e., Axiom K 3 is also true and $Q$ is a classical probability function.

Besides, $P$ is represented by $Q$ with the homomorphism $f: [0, m] \to [0, 1]$ defined as $f(n) = n/k$.

Theorem is proved.

It is necessary to remark that in spite of the possibility to faithfully represent any inflated probability function by a classical probability function, utilization of inflated probability better correlates with informal probabilistic reasoning and makes such reasoning more efficient (Burgin and Meissner, 2012).

In some probability structures, order is generated by probabilities of elementary events.

**Definition 3.5.** A positive probability structure $P: \mathbf{F} \to X$ is called *elementary*, if for any elements $\{u_1, u_2, u_3, \ldots, u_q\}$ and $\{v_1, v_2, v_3, \ldots, v_t\}$ from $\mathbf{F}$ (here all $u_1, u_2, u_3, \ldots, u_q$ and $v_1, v_2, v_3, \ldots, v_t$ belong to $U$), we have

$$P(\{u_1, u_2, u_3, \ldots, u_q\}) \succ P(\{v_1, v_2, v_3, \ldots, v_t\})$$

if and only if $q = t$ and for any $1 \leq i \leq t$ either $u_i = v_i$ or $P(u_i) \succ P(v_i)$ and the latter is true for, at least, one $i$.

**Theorem 3.2.** It is possible to faithfully represent any total elementary positive probability structure $P$ by a classical probability function.

*Proof.* Let us consider a total positive probability structure $P: 2^U \to X$. It is possible to extend the strict partial order $\succ$ in $X$ to a total (linear) order in $X$, which we denote by the same symbol order $\succ$ in $X$. Indeed, the strict partial order $\succ$ in $X$ corresponds to a directed acyclic graph $G = (V, E)$ in which vertices are labeled by elements of $X$. Any directed acyclic graph has the topological ordering of its vertices, which defines a total (linear) order in X (Kleinberg and Tardos, 2006). In such a way, X becomes totally ordered.

In addition we see that the set U is equal to the union $\bigcup_{i=0}^{r} U_i$ where $P(x) = P(y)$ for any two elements (elementary events) x and y from the same set $U_i$ and $P(x) \neq P(z)$ for any two elements (elementary events) x and z from different sets $U_i$. Besides, $U_0$ consists of all elements z from U such that $P(z) = P(\emptyset)$ and we enumerate sets $U_i$ so that $i < j$ when $P(x) < P(z)$ for elements (elementary events) x from $U_i$ and z from $U_{ij}$. For an element (elementary event) z from the set $U_0$, we define its probability $p(x) = 0$.

As we consider only finite sets of events, all sets $U_i$ are finite and let the number of elements in $U_i$ be $|U_i| = n_i$. Then we take the sum

$$n_1 + 2n_2 + 3n_3 + \ldots + rn_r = n \qquad (2)$$

and for an element (elementary event) x from the set $U_i$ when $i > 0$, we define its probability as

$$p(x) = i/n$$

Then we extend the mapping $p: U \to X$ to the function $P: 2^U \to X$ by the following rule:

If $x \in 2^U$, then $x = \{ u_1, u_2, u_3, \ldots, u_k \}$ where all $u_k$ belong to U and

$$p(x) = p(u_1) + p(u_2) + p(u_3) + \ldots + p(u_k)$$

This defines the mapping $f: X \to \mathbf{R}$ and allows us to show that the function $p: 2^U \to \mathbf{R}$ is a classical probability function. To show this, we have to check Axioms K1 – K3.

K 1. Any event is a sum of elementary events. Thus, $p(A) \geq 0$, for all $A \in 2^X$ because for any elementary event x from the set U, $p(x) = i/n \geq 0$ when $x \in U_i$ with $i > 0$ and $p(x) = 0$ when $x \in U_0$.

K 2. As $U = \{w_1, w_2, w_3, \ldots, w_n\}$, we have

$p(\Omega) = p(w_1) + p(w_2) + p(w_3) + \ldots + p(w_n) = n_0(0/n) + n_1(1/n) + n_2(2/n) + n_3(3/n) + \ldots + n_r(r/n) = 1$

K 3. Let us consider two disjoint subsets A and B from U, i.e., $A \cap B = \emptyset$. Then

$$A = \{u_1, u_2, u_3, \ldots, u_q\}$$
$$B = \{v_1, v_2, v_3, \ldots, v_t\}$$

and

$$A \cup B = \{u_1, u_2, u_3, \ldots, u_q, v_1, v_2, v_3, \ldots, v_t\}$$

By definition, we have

$P(A \cup B) = p(u_1) + p(u_2) + p(u_3) + \ldots + p(u_n) + p(v_1) + p(v_2) + p(v_3) + \ldots + p(v_n) = P(A) + P(B)$

i.e., $P(A \cup B) = P(A) + P(B)$.

By construction, the mapping $f: X \to \mathbf{R}$ is a homomorphism of the partially ordered set $X$ into the ordered set $\mathbf{R}$ because when $i < j$, we $p(x) == i/n < p(z) = j/n$ for elements (elementary events) $x$ from $U_i$ and $z$ from $U_{ij}$.

Theorem is proved.

.**Remark 3.1.** The representation of positive probability structures described in Theorem 3.2 is not unique. Indeed, it is also possible to use the following representation taking the sum

$$2n_1 + 4n_2 + 6n_3 + \ldots + 2rn_r = m \qquad (3)$$

and for an element (elementary event) $x$ from the set $U_i$ ($= 0, 1, \ldots ,$ ), defining its probability as

$$p(x) = 2i/n$$

Then we extend the mapping $p: U \to X$ to the function $P: \mathbf{F} \to X$ by the following rule:

If $x \in \mathbf{F}$, then $x = \{ u_1, u_2, u_3, \ldots, u_k \}$ where all $u_k$ belong to $U$ and

$$p(x) = p(u_1) + p(u_2) + p(u_3) + \ldots + p(u_k)$$

**Remark 3.2.** Results from (Kraft, et al, 1959) show that there are probability structures for which it is impossible to find a representation by classical probability.

.**Remark 3.3.** It is necessary to understand that reduction of one structure to another one does not eliminate necessity in the first structure. Usually, in some situations, the first structure might be more relevant and/or more efficient, while in other situations, the second structure can provide more possibilities to deal with the problem under consideration.

## 4. Conclusion

We developed the theory of qualitative probability based on probability structures obtaining various properties of these structures. We demonstrated that classical probability and inflated probability are positive probability structures. In particular, we explored representations of some probability structures by other probability structures, for example, inflated probability by classical probability or some classes of total probability structures by classical probability.

Obtained results allow us to formulate open problems and interesting directions for further research.

It would be appealing to study independence of axioms considered in this paper.

One more worthy of note question is whether the following axiom is valid in any positive probability structure.

**Positive Complementarity Axiom (PCA)**: If $P(A) \succ P(B)$, then it is not $P(CA) \succ P(CB)$ where $CX$ is the complement of the event $X$.

It would be interesting and important to study conditional probability structures.

Qualitative probabilities in the form of comparative probabilities have been studied in the context of positive probability. The same is true for probability structures studied in this paper. However, intensive utilization and exploration of negative probabilities poses the problem of introduction and investigation of qualitative probabilities that include negative probabilities.

In this paper, we studied only comparative axioms for qualitative probabilities. At the same time, researchers also considered axioms for independence of events and axioms for uncertainty of experiments. Thus, an essential problem is to study these axioms in the context of probability structures.